\documentclass[12pt]{elsarticle}
\usepackage{amsfonts}
\usepackage{amsmath}
\usepackage{graphicx}
\usepackage{epsfig}
\usepackage{amssymb}
\usepackage{amstext}
\usepackage{tikz}
\newtheorem{theorem}{Theorem}[section]
\newtheorem{lemma}{Lemma}[section]

\newtheorem{corollary}{Corollary}[section]

\newenvironment{PrfFact}{{\bf Proof }}{{\hfill\tiny{$\blacksquare$\\}}}

\begin{document}
\begin{frontmatter}

\title{On Barnette's Conjecture and $H^{+-}$ property}

\author{Jan~Florek}
\ead{jan.florek@ue.wroc.pl}

\address{Institute of Mathematics,\\
University of Economics \\
ul. Komandorska 118/120
\\ 53--345 Wroc{\l}aw, Poland}
\begin{abstract}
A conjecture of Barnette  states that every $3$-connected cubic bipartite plane graph has a Hamilton cycle, which is equivalent to the statement that every simple even plane triangulation admits a partition of its vertex set into two subsets so that each induces a tree.

Let $G$ be a simple even plane triangulation and suppose that ${V_1, V_2, V_3}$ is a $3$-coloring of the vertex set of $G$. Let $B_{i}$, $i = 1, 2, 3$, be the set of all vertices in $V_i$ of the degree at least $6$. We prove that if induced graphs $G[B_1 \cup B_2]$ and $G[B_1 \cup B_3]$ are acyclic, then the following properties are satisfied:
\\[6pt]
(1) For every path $abc$ there is possible to partition the vertex set of $G$ into two subsets so that each induces a tree, and one of them contains the edge $ab$ and avoids the vertex $c$,
\\[6pt]
(2)  For every path $abc$  with vertices $a$, $c$ of the same color there is possible to partition the vertex set of $G$ into two subsets so that each induces a tree, and one of them contains the path $abc$.
 \end{abstract}

\end{frontmatter}
\section{Introduction}

All graphs considered in this paper are simple plane graphs. We use \cite{flobar1} for general reference.

Let ${\cal P}$ be the family of all 3-connected cubic bipartite plane graphs. Barnette (\cite{flobar11}, Problem~5) conjectured that every graph in ${\cal P}$ is Hamiltonian. In \cite{flobar4}, Goodey proved that if a graph in ${\cal P}$ has only faces with 4 or 6 sides, then it is Hamiltonian. Feder and Subi \cite{flobar2} generalized Goodey`s result by showing that when the faces of a graph in ${\cal P}$ are 3-colored, if two of the three color classes contain only faces with either 4 or 6 sides the conjecture holds. Holton, Manvel and McKay \cite{flobar6} used computer search to confirm Barnette's conjecture for graphs up to 64 vertices, and also to confirm the following properties $H^{+-}$ and $H^{--}$ for graphs up to 40 vertices:
\\[6pt]
$H^{+-}$ \quad If any two edges are chosen on the same face, there is a Hamilton cycle through one and avoiding the other,
\\[6pt]
$H^{--}$ \quad If any two edges are chosen which are an even distance apart on the same face, there is a Hamilton cycle which avoids both.
\\[6pt]
Kelmans \cite{flobar7}  proved that Barnette's conjecture holds if and only if every graph in ${\cal P}$ has the property $H^{+-}$.

Let ${\cal E}$ be the dual family of ${\cal P}$. Thus ${\cal E}$ is the family of of all simple even plane triangulations. It is known that  Barnette's conjecture is equivalent to the following statement : every graph in ${\cal E}$ has a partition of its vertex set into two subsets so that each induces a tree (see Stein \cite{flobar10}). The dual approach to Barnette's conjecture was used by the author \cite{flobar3} and Lu \cite{flobar8}, \cite{flobar9}. Let $G \in {\cal E}$, and suppose that ${V_1, V_2, V_3}$ is a $3$-coloring of $V(G)$. Notice that properties $H^{+-}$ and $H^{--}$ are equivalent to the following dual versions:
\\[6pt]
(1) For every path $abc$ there is possible to partition the vertex set of $G$ into two subsets so that each induces a tree, and one of them contains the edge $ab$ and avoids the vertex $c$,
\\[6pt]
(2)  For every path $abc$  with vertices $a$, $c$ of the same color there is possible to partition the vertex set of $G$ into two subsets so that each induces a tree, and one of them contains the path $abc$.
\\[6pt]
Denote by $B_{i}$, $i = 1, 2, 3$, the set of all vertices in $V_i$ of the degree at least $6$ in $G$. The purpose of this paper is to prove Theorem~\ref{theorem2.1}: if induced graphs $G[B_1 \cup B_2]$ and $G[B_1 \cup B_3]$ are acyclic, then the properties (1)--(2) are satisfied (which generalize results obtained by the author in \cite{flobar3}). Notice that if one of the above induced graphs is acyclic, then by  Corollary~\ref{corollary2.1} the dual graph $G^*$ is Hamiltonian.



\section{Main results}
Let $G$ be a simple even plane triangulation, and suppose that ${V_1, V_2, V_3}$ is a $3$-coloring of $V(G)$; it is well know that such coloring exists (see Heawood  \cite{flobar5}). We say that a vertex $v$ is \textsl{small} (\textsl{big}) in $G$ if $deg_{G}(v) = 4$ ($deg_{G}(v) \geq 6$, respectively). Denote by $B_i$ (or $S_i$) the set of all big (or small, respectively) vertices in $V_i$, $i = 1, 2, 3$.
\begin{lemma}\label{lemma2.1}
Let $X$, $Y$ be disjoint sets of vertices in the graph $G$ such that:

(a) $B_1 \subseteq X$, $B_2 \subseteq Y$, and $B_3 \subseteq X \cup Y$,

(b) for every induced path with big ends and small its all inner vertices, the set of inner vertices is contained in or is disjoint with $X \cup Y$.

If the induced graphs $G[X]$ and $G[Y]$ are acyclic, then it is possible to partition the vertex set of $G$ into two subsets $S$, $T$ so that each induces a tree, $X \subseteq S$, and $Y \subseteq T$.
\end{lemma}
\begin{PrfFact}
If $G$ has at most two big vertices, then theorem yields. Therefore we assume that $G$ has at least three big vertices. Let $v$ be a small vertex which not belongs to $X \cup Y$. There are two possibilities:

(i) $v$ belongs to a maximum induced path $P$, with big ends in $G$ (say $a$, $b$), and small its all inner vertices. Moreover, there exist two big vertices (say $c$, $d$) which are adjacent to all vertices of the path $P$.

(ii) $v$ belongs to a subgraph $H \cong K_4$ (octahedron).

Case (i).  Notice that vertices $c$ and $d$ belong to the same color class. Hence, by (a), if $c \in X$ and $d \in Y$, then $c,d \in V_3$. Thus the vertices of $P$ belong to $V_1 \cup V_2$.  By (b), $X'= X \cup (V(P) \cap S_1)$ and $Y' = Y \cup (V(P)\cap S_2)$ are disjoint. Notice that for $X'$ and $Y'$ the conditions (a)--(b) are satisfied, and the induced graphs $G[X']$ and $G[Y']$ are disjoint and acyclic.

If $c,d \in X$ and doesn't exist a $cd$-path consisting of  vertices in $X$, we choose a small vertex $e$ belonging to the path $P$. Then for $X'= X \cup \{e\}$ and $Y' = Y \cup V(P - e)$ the conditions (a)--(b) are satisfied, and the induced graphs $G[X']$ and $G[Y']$ are disjoint and acyclic.

If $c,d \in X$ and there is a $cd$-path consisting of  vertices in $X$, then for $X' = X$ and $Y' = Y \cup (V(P)\backslash X)$ the conditions (a)--(b) are satisfied, and the induced graphs $G[X']$ and $G[Y']$ are disjoint and acyclic.

Case (ii). The subgraph $H$ is bounded by a $3$-cycle (say $C$) consisting of big vertices in $G$, and has three small vertices in $G$ (say $x$, $y$ and $z$). Assume that $C$ has two vertices in $Y$ which both are adjacent to $x$. Then for $X' = X \cup \{x,z\}$ and $Y' = Y \cup \{y\}$ the conditions (a)--(b) are satisfied, and the induced graphs $G[X']$ and $G[Y']$ are disjoint and acyclic.

Continuing the process for the sets $X'$, $Y'$ as for the sets $X$, $Y$ we obtain a partitioning of $V(G)$ into two subsets $S$, $T$ each of which induces a tree,  $X \subseteq X' \subseteq S$, and $Y \subseteq Y' \subseteq T$.
\end{PrfFact}
\begin{corollary}\label{corollary2.1}
Let $X$, $Y$ be disjoint sets of big vertices in $G$ such that
$B_1 \subseteq X$, $B_2 \subseteq Y$, and $B_3 \subseteq X \cup Y$.
If the induced graphs $G[X]$ and $G[Y]$ are acyclic, then it is possible to partition the vertex set of $G$ into two subsets $S$, $T$ so that each induces a tree.
\end{corollary}
If $G$ is $4$-connected and has at least three big vertices, then for every small  vertex $s$, which is adjacent to a small vertex in $G$, there exists just one maximum induced path $P_{s}$ containing $s$, with big ends and small its all inner vertices. If a small vertex $s$ is adjacent with four big vertices in $G$, then there are two such paths of the length $2$ containing the vertex $s$. In such case we indicate which one of these two paths has been chosen as $P_s$.

\begin{lemma}\label{lemma2.2}
Let $G$ be 4-connected with at least three big vertices, and suppose that the induced graphs $G[B_1 \cup B_2]$ and $G[B_1 \cup B_3]$ are acyclic. For every path $abc$ with vertices $a$, $c$ of the same color there are disjoint sets $X ,Y \subset V(G)$ such that :

(a) the induced graphs $G[X]$ and $G[Y]$ are acyclic, and one of them contains the edge $ab$ and avoids the vertex $c$,

(b) there is a permutation $(i,j,k)$ of $\{1,2,3\}$ and a vertex $v \in B_j$ such that $B_i \cup B_j \backslash \{v\} \subseteq X$ and $B_k \cup\{v\} \subseteq Y$, or $B_i \cup B_j \subseteq X$ and $B_k \subseteq Y$,

(c) for every induced path with big ends in $G$ and small its all inner vertices, the set of inner vertices is contained in or is disjoint with $X \cup Y$.
\end{lemma}
\begin{PrfFact}
Assume that the induced graphs $G[B_1 \cup B_2]$ and $G[B_1 \cup B_3]$ are acyclic. Fix a path $abc$. Suppose first that the vertices of the path $abc$ are big. By symmetry of the assumption, it is sufficient to consider the following cases:
\\
($\alpha$) $a \in B_1$, $b \in B_2$ and $c \in B_1$,
\\
($\beta$) $a \in B_2$, $b \in B_1$ and $c \in B_2$,
\\
($\gamma$) $a \in B_2$, $b \in B_3$ and $c \in B_2$.
\\
We set

$X = (B_1 \cup B_2) \backslash \{c\}$ and $Y = B_3 \cup \{c\}$, for the case ($\alpha$) and ($\beta$),

$X = B_1 \cup B_2 \backslash \{a\}$ and $Y = B_3 \cup \{a\}$, for the case ($\gamma$).
\\
Since $G[B_1 \cup B_2]$ is acyclic and $B_3$ is independent, $G[X]$ and $G[Y]$ are acyclic too. Notice that $G[X]$ ($G[Y]$) contains the edge $ab$ and avoids the vertex $c$ for the case  ($\alpha$) and ($\beta$) (for the case ($\gamma$), respectively). Hence, the conditions (a)--(c) are satisfied.

Let now the path $abc$ has a small vertex in $G$. By symmetry of the assumption, it is sufficient to consider the following cases :
\\
($\alpha, 1$) $a \in S_1$, $b \in B_2$, $b \notin V(P_a)$ and $c \in S_1$,
\\
($\alpha, 2$) $a \in S_1$, $b \in B_2$, $b \in V(P_a)$ and $c \in S_1$,
\\
($\alpha, 3$) $a \in V_1$, $b \in S_2$ and $c \in V_1$,
\\
($\alpha, 4$) $a \in B_1$, $b \in B_2$ and $c \in S_1$,
\\
($\alpha, 5$) $a \in S_1$, $b \in B_2$ and $c \in B_1$,
\\
($\beta, 1$) $a \in S_2$, $b \in B_1$ and $c \in B_2$,
\\
($\beta, 2$) $a \in B_2$, or $b \in S_1$, or $c \in S_2$,
\\
($\gamma, 1$) $a \in S_2$, $b \in B_3$ and $c \in S_2$,
\\
($\gamma, 2$) $a \in V_2$, $b \in S_3$ and $c \in V_2$,
\\
($\gamma, 3$) $a \in B_2$, $b \in B_3$ and $c \in S_2$,
\\
($\gamma, 4$) $a \in S_3$, $b \in B_2$ and $c \in B_3$.

Case ($\alpha,1$). There is a vertex $d \in B_2 \backslash \{b\}$ which is adjacent to every vertex of the path $P_a$. Since $B_1 \cup B_2$ is acyclic there exists at most one $bd$-path consisting of vertices in $B_1 \cup B_2$. If such path exists, let $be$, $e \in B_1$, be an edge of this path (if such path not exist we replace $\{e\}$ by the empty set). If $b \notin V(P_c)$ and $P_a \neq P_c$, then we set
$$X = B_1 \backslash \{e\} \cup B_2 \cup \{a\}, \quad
Y = B_3 \cup  V(P_a - a) \cup (V(P_c)\backslash B_1) \cup \{e\} \quad \hbox{and} \quad v:= e.$$
Since $G$ is $4$-connected, $e$ is not adjacent to both ends of $P_c$. Hence, $G[X]$ and $G[Y]$ are acyclic, because $G[B_1 \cup B_2]$ is acyclic and $B_3$ is independent. Notice that $X$, $Y$ are disjoint, $ab \subset G[X]$, ${(B_1 \backslash \{e\}) \cup B_2} \subset X$, $B_3 \cup \{e\} \subset Y$, and $X \cup Y = B_1 \cup B_2 \cup B_3 \cup V(P_a) \cup V(P_c)$. Hence, conditions (a)--(c) yield.

If $b \in V(P_c)$, then we set
$$X = B_1 \backslash \{e\} \cup B_2 \cup \{a\} \cup V(P_c - c), \quad
Y = B_3 \cup  V(P_a - a) \cup  \{c, e\} \quad \hbox{and} \quad v:= e.$$
Since $G$ is $4$-connected, $e$ is not adjacent to both vertices in $B_3$ which are adjacent to every vertex of $P_c$. Hence, $G[X]$ and $G[Y]$ are acyclic and conditions (a)--(c) yield.

If $P_a = P_c$, then we set
$$X = B_1 \backslash \{e\} \cup B_2 \cup \{a\}, \quad
Y = B_3 \cup  V(P_a - a) \cup \{e\} \quad \hbox{and} \quad v:= e.$$

Case ($\alpha,2$).
If $a$ is adjacent with four big vertices in $G$, then  the case ($\alpha,2$) follows from the case ($\alpha,1$). Thus we assume that the path $P_a$ has the length at least $3$. Hence, there exists a vertex $d \in V(P_a) \cap S_2$. If $b \notin V(P_c)$, then we set
$$X = B_1 \cup B_2 \cup V(P_a - d) \quad \hbox{and} \quad Y = B_3 \cup \{d\} \cup V(P_c)\backslash B_1.$$
If $b \in V(P_c)$, then we set
$$X = B_1 \cup B_2 \cup V(P_a - d) \cup  V(P_c - c) \quad \hbox{and} \quad Y = B_3 \cup  \{c, d\}.$$

Case ($\alpha,3$). If  $a,c \notin V(P_b)$, then we set
$$X = B_1 \backslash \{c\} \cup B_2 \cup (V(P_b) \cap S_2), \quad Y = B_3 \cup \{c\} \cup (V(P_b) \cap S_3) \quad \hbox{and} \quad v:= c.$$
Let $a,c \in V(P_b)$. If $c$ is small (big) in $G$,  then we set
$$X = B_2 \cup (B_1 \cup V(P_b)) \backslash \{c\}, \quad Y = B_3 \cup \{c\}, \quad (\hbox{and } v:= c,  \hbox{ respectively}).$$

Case ($\alpha,4$). If $b \notin V(P_c)$, then we set
$$X = B_1 \cup B_2  \quad \hbox{and} \quad Y = B_3 \cup  V(P_c) \backslash B_1.$$
If $b \in V(P_c)$, then we set
$$X = B_1 \cup B_2\cup V(P_c - c) \quad \hbox{and} \quad Y = B_3 \cup \{c\}.$$

Case ($\alpha, 5$). Suppose that $b \notin V(P_a)$. Then we set
$$X = B_1 \cup B_3 \cup V(P_a - a) \quad \hbox{and} \quad Y = B_2 \cup  \{a\}.$$
Suppose that $b \in V(P_a)$. Then we set
$$X = B_1 \cup B_3 \quad \hbox{and} \quad Y = B_2 \cup V(P_a) \backslash B_1.$$

Case ($\beta$,1).
Suppose that $b \notin V(P_a)$. Since $G[B_1 \cup B_2]$ and $G[B_1 \cup B_3]$ are acyclic, exactly one end of the path $P_a$ (say $d$) belongs to $B_3$. We set
$$X = B_1 \cup B_3 \backslash \{d\} \cup \{a\}, \quad Y = B_2 \cup V(P_a - a) \quad \hbox{and} \quad v:= d.$$
There is a vertex $e \in B_1\backslash \{b\}$ which is adjacent to every vertex of $P_a$. Since $G[B_1 \cup B_3]$ is acyclic, $bde$ is the only $be$-path consisting of vertices in $B_1 \cup B_3$. Hence, $G[X]$ and $G[Y]$ are acyclic, because $G[B_1 \cup B_3]$ is acyclic and $B_2$ is independent. Notice that $X$, $Y$ are disjoint, $ab \subset G[X]$, ${B_1 \cup B_3} \backslash \{d\} \subset X$, $B_2 \cup \{d\} \subset Y$, and $X \cup Y = B_1 \cup B_2 \cup B_3 \cup V(P_a)$. Therefore the conditions (a)--(c) yield.

Suppose that $b \in V(P_a)$. Since $G[B_1 \cup B_3]$ is acyclic, there is a vertex $d \in V(P_a) \cap S_1$. Since $G[B_1 \cup B_2]$ is acyclic there exists at most one path connecting the ends of the path $P_a$ and consisting of vertices in $B_1 \cup B_2$. If such path exists, we set
$$X = B_1 \cup B_2 \backslash \{c\} \cup V(P_a - d), \quad
Y = B_3 \cup  \{c, d\} \quad \hbox{and} \quad v:= c.$$
If such path not exists, we set
$$X = B_1 \cup B_2  \backslash \{c\} \cup V(P_a), \quad
Y = B_3 \cup  \{c\} \quad \hbox{and} \quad v:= c.$$

Case ($\beta$,2). For the cases ($\alpha,1$)--($\alpha,4$), Lemma~\ref{lemma2.2} follows from the assumption that the induced graph $G[B_1 \cup B_2]$ is acyclic. Hence,  for the case ($\beta,2$), Lemma~\ref{lemma2.2} yields.

Case ($\gamma$,1). Suppose that $b \in V(P_a) \cap V(P_c)$. There is a vertex $d\in B_1$ which is adjacent to every vertex of $P_c$. Then we set
$$X = B_1 \backslash \{d\} \cup B_2 \cup (V(P_c) \cap S_2),$$
$$Y = B_3\cup \{d\} \cup (V(P_a) \backslash B_2) \cup (V(P_c) \cap S_3) \quad \hbox{and} \quad v:= d.$$
Since $G[B_1 \cup B_3]$ is acyclic, the path $P_a$  has one end belonging to $B_2$. Hence, $G[X]$ and $G[Y]$ are acyclic, because $G[B_1 \cup B_2]$ is acyclic and $B_3$ is independent. Notice that $X$, $Y$ are disjoint, $ab \subset G[Y]$, $B_1 \backslash \{d\} \cup B_2 \subset X$, $B_3 \cup \{d\} \subset Y$, and $X \cup Y = B_1 \cup B_2 \cup B_3 \cup V(P_a) \cup V(P_c)$. Hence, conditions (a)--(c) yield.

Let $b \notin V(P_a)$ and $b \in V(P_c)$. Since $G[B_1 \cup B_3]$ is acyclic, the path $P_c$  has one end (say $d$) belonging to $B_2$. Then we set
$$X = B_1 \cup B_2  \backslash \{d\} \cup V(P_a - a) \cup \{c\}, \quad Y = B_3 \cup \{a\} \cup V(P_c - c) \quad \hbox{and} \quad v:= d.$$
There are vertices $e,f \in B_1$ which are adjacent to every vertex of $P_c$. Since $G[B_1 \cup B_2]$ is acyclic, $edf$ is the only $ef$-path consisting of vertices in $B_1 \cup B_2$. Hence, $G[X]$ and $G[Y]$ are acyclic, because $G[B_1 \cup B_2]$ is acyclic and $B_3$ is independent. Notice that $X$, $Y$ are disjoint, $ab \subset G[Y]$, ${B_1 \cup B_2} \backslash \{d\} \subset X$, $B_3 \cup \{d\} \subset Y$, and $X \cup Y = B_1 \cup B_2 \cup B_3 \cup V(P_a) \cup V(P_b)$. Therefore the conditions (a)--(c) yield.

Let $b \in V(P_a)$ and $b \notin V(P_c)$. Since $G[B_1 \cup B_3]$ is acyclic,
there is a vertex $d \in V(P_c) \cap S_1$. Then we set
$$X = B_1 \cup B_2 \cup V(P_c - d) \quad \hbox{and} \quad Y = B_3 \cup \{d\} \cup (V(P_a)  \backslash B_2).$$
Let $b \notin V(P_a) \cup V(P_c)$ and $P_a \neq P_c$.  Since $G[B_1 \cup B_3]$ is acyclic, there is a vertex $d \in {V(P_c) \cap S_1}$. Then we set
$$X = B_1 \cup B_2 \cup V(P_a - a) \cup V(P_c -d) \quad \hbox{and} \quad Y = B_3 \cup \{a, d\}.$$
If $b \notin V(P_a)$ and  $P_a = P_c$, then we set
$$X = B_1 \cup B_2 \cup V(P_a - a) \quad \hbox{and} \quad Y = B_3 \cup \{a\}.$$

Case ($\gamma$,2).
If $a,c \notin V(P_b)$, then we set
$$X = B_1 \cup B_2 \backslash \{a\} \cup (V(P_b) \cap S_1), \quad Y = B_3 \cup \{a\} \cup (V(P_b) \cap S_3) \quad \hbox{and} \quad v:= a.$$
Let $a,c \in V(P_b)$.  Since $G[B_1 \cup B_2]$ and  $G[B_1 \cup B_3]$, exactly one end of the path $P_b$ (say $d$) belongs to $B_2$. If $d \neq c$, then we set
$$X = B_1 \cup B_2 \backslash \{d\} \cup \{c\}, \quad Y = B_3 \cup V(P_b - c) \quad \hbox{and} \quad v:= d.$$
There are vertices $e,f \in B_1$ which are adjacent to every vertex of $P_b$. Since $G[B_1 \cup B_2]$ is acyclic, $edf$ is the only $ef$-path consisting of vertices in $B_1 \cup B_2$. Hence, $G[X]$ and $G[Y]$ are acyclic, because $G[B_1 \cup B_2]$ is acyclic and $B_3$ is independent. Notice that $X$, $Y$ are disjoint, $ab \subset G[Y]$, $B_1 \cup B_2 \backslash \{d\} \subset X$, $B_3 \cup \{d\} \subset Y$, and $X \cup Y = B_1 \cup B_2 \cup B_3 \cup V(P_b)$. Hence, conditions (a)--(c) yield.
If $d = c$, then we set
$X = B_1 \cup B_2$ and $Y = B_3 \cup V(P_b - c)$.

Case ($\gamma$,3). Suppose that $b \in V(P_c)$. Then we set
$$X = B_1 \cup B_3 \backslash \{b\} \cup  \{c\}, \quad Y = B_2 \cup V(P_c - c) \quad \hbox{and} \quad v:= b.$$
There are vertices $e,f \in B_1$ which are adjacent to every vertex of $P_c$. Since $G[B_1 \cup B_3]$ is acyclic, $ebf$ is the only $ef$-path consisting of vertices in $B_1 \cup B_3$. Hence, $G[X]$ and $G[Y]$ are acyclic, because $G[B_1 \cup B_3]$ is acyclic and $B_2$ is independent. Notice that $X$, $Y$ are disjoint, $ab \subset G[Y]$, $B_1 \cup B_3 \backslash \{b\} \subset X$, $B_2 \cup \{b\} \subset Y$, and $X \cup Y = B_1 \cup B_2 \cup B_3 \cup V(P_c)$. Hence, conditions (a)--(c) yield.

Suppose that $b \notin V(P_c)$. Since $G[B_1 \cup B_3]$ is acyclic, there is a vertex $d \in V(P_c) \cap S_1$. Since $G[B_1 \cup B_2]$ is acyclic there exists at most one path connecting ends of the path $P_c$ and consisting of vertices in $B_1 \cup B_2$. If such path exists, we set
$$X = B_1 \cup B_2 \backslash \{a\} \cup V(P_c - d), \quad
Y = B_3 \cup  \{a, d\} \quad \hbox{and} \quad v:= a.$$
If such path not exists, we set $X = B_1 \cup B_2  \backslash \{a\} \cup V(P_c)$ and $Y = B_3 \cup  \{a\}$.

Case ($\gamma$,4). For the case ($\alpha,5$), Lemma~\ref{lemma2.3} follows from the assumption that the induced graph $G[B_1 \cup B_3]$ is acyclic. Hence, for the case ($\gamma,4$), Lemma~\ref{lemma2.3} yields.
\end{PrfFact}
\begin{lemma}\label{lemma2.3}
Let $G$ be 4-connected with at least three big vertices, and suppose that the induced graphs $G[B_1 \cup B_2]$ and $G[B_1 \cup B_3]$ are acyclic. For every path $abc$ with vertices $a$, $c$ of the same color there are disjoint sets $X ,Y \subset V(G)$ such that :

(a) the induced graphs $G[X]$ and $G[Y]$ are acyclic and one of them contains the path $abc$,

(b) there is a permutation $(i,j,k)$ of $\{1,2,3\}$ and a vertex $v \in B_j$ such that $B_i \cup B_j \backslash \{v\} \subseteq X$ and $B_k \cup\{v\} \subseteq Y$, or $B_i \cup B_j \subseteq X$ and $B_k \subseteq Y$,

(c) for every induced path with big ends in $G$ and small its all inner vertices, the set of inner vertices is contained in or is disjoint with $X \cup Y$.
\end{lemma}
\begin{PrfFact}  Assume that the induced graphs $G[B_1 \cup B_2]$ and $G[B_1 \cup B_3]$ are acyclic. Fix a path $abc$. Suppose first that  vertices of the path $abc$ are big in $G$ . By symmetry of the assumption, it is sufficient to consider the following cases:
\\
($\alpha$) $a \in B_1$, $b \in B_3$ and $c \in B_1$,
\\
($\beta$) $a \in B_3$, $b \in B_1$ and $c \in B_3$,
\\
($\gamma$) $a \in B_3$, $b \in B_2$ and $c \in B_3$.
\\
We set

$X = B_1 \cup B_3$ and $Y = B_2$, for the cases ($\alpha$) and ($\beta$),

$X = B_1 \cup B_2 \backslash \{b\}$ and $Y = B_3 \cup \{b\}$, for the case ($\gamma$).
\\
Since $G[B_1 \cup B_2]$ and $G[B_1 \cup B_3]$ are acyclic and $B_i$ is independent, $G[X]$ and $G[Y]$ are acyclic too. Notice that the path $abc$ is contained in $G[X]$ ($G[Y]$) for the cases  ($\alpha$) and ($\beta$) (for the case ($\gamma$), respectively). Hence, the conditions (a)--(c) are satisfied.

Let now the path $abc$ has a small vertex in $G$. By symmetry of the assumption, it is sufficient to consider the following cases :
\\
($\alpha, 1$) $a \in S_1$, $b \in B_3$ and $c \in S_1$,
\\
($\alpha, 2$) $a \in V_1$, $b \in S_3$ and $c \in V_1$,
\\
($\alpha, 3$) $a \in S_1$, $b \in B_3$ and $c \in B_1$,
\\
($\beta, 1$) $a \in S_3$, $b \in B_1$ and $c \in S_3$,
\\
($\beta, 2$) $a \in V_3$, $b \in S_1$ and $c \in V_3$,
\\
($\beta, 3$) $a \in S_3$, $b \in B_1$ and $c \in B_3$,
\\
($\gamma$) $a \in V_3$, $b \in V_2$ and $c \in V_3$.

Case ($\alpha, 1$). If $b \notin V(P_a) \cup V(P_c)$ and $P_a \neq P_c$, then we set
$$X = B_1 \cup B_2 \cup V(P_a - a) \cup V(P_c -c) \quad \hbox{and} \quad Y = B_3 \cup \{a, c\}.$$
If $b \in V(P_a) \cap V(P_c)$, then we set
$$X = B_1 \cup B_2 \quad \hbox{and} \quad Y = B_3 \cup (V(P_a) \cup V(P_c))\backslash B_1.$$
If $b \notin V(P_a)$ and $b \in V(P_c)$, then we set
$$X = B_1 \cup B_2 \cup V(P_a - a)\quad \hbox{and} \quad Y = B_3 \cup \{a\} \cup V(P_c)\backslash B_1.$$
Suppose that $b \notin V(P_a)$ and  $P_a = P_c$.
There exists a vertex $d \in B_3 \backslash \{b\}$ which is adjacent to every vertex of the path $P_a$. Then we set
$$X = B_1 \cup B_3\backslash \{d\} \cup (V(P_a) \cap S_1), \quad Y = B_2 \cup \{d\} \cup (V(P_a) \cap S_2) \quad \hbox{and} \quad v := d.$$

Case ($\alpha, 2$). Let  $a,c \notin V(P_b)$.
Since $G[B_1 \cup B_2]$ and $G[B_1 \cup B_3]$ are acyclic, exactly one end of the path $P_b$ (say $d$) belongs to $B_2$. Then we set
$$X = B_1 \cup B_2 \backslash \{d\} \cup \{b\}, \quad Y = B_3 \cup V(P_b - b)\quad \hbox{and} \quad v:= d.$$
Since $G[B_1 \cup B_2]$ is acyclic, $adc$ is the only $ac$-path with vertices in $B_1 \cup B_2$. Hence, $G[X]$ and $G[Y]$ are acyclic, because $G[B_1 \cup B_2]$ is acyclic and $B_3$ is independent. Notice that $X$, $Y$ are disjoint, $abc \subset G[X]$, ${B_1 \cup B_2} \backslash \{d\} \subset X$, $B_3 \cup \{d\} \subset Y$, and $X \cup Y = B_1 \cup B_2 \cup B_3 \cup V(P_b)$. Therefore the conditions (a)--(c) yield.

Suppose that $a,c \in V(P_b)$. Since $G[B_1 \cup B_2]$ is acyclic, at least one end of the path $P_b$ (say $d$) belongs to $B_3$. Then we set
$$X = B_1 \cup (B_3 \cup V(P_b)) \backslash \{d\}, \quad Y = B_2 \cup \{d\} \quad \hbox{and} \quad v:= d.$$

Case ($\alpha, 3$). If $b \notin V(P_a)$, then we set
$$X = B_1 \backslash \{c\} \cup B_2 \cup V(P_a - a), \quad Y = B_3 \cup  \{a, c\} \quad \hbox{and} \quad v:= c.$$
Suppose that  $b \in V(P_a)$.  Let $d \neq b$ be the end of the path $P_a$. Notice that $d \in B_1 \cup B_3$. Then we set
$$X = (B_1 \cup B_3 \cup V(P_a)) \backslash \{d\}, \quad Y = B_2 \cup \{d\} \quad \hbox{and} \quad v: = d.$$

Case ($\beta, 1$). If $b \notin V(P_a)$ and $b \in V(P_c)$, then we set
$$X = B_1 \backslash \{b\} \cup B_2 \cup (V(P_a) \cap S_2), \quad Y = B_3 \cup (V(P_a) \cap S_3) \cup V(P_c) \quad \hbox{and} \quad v:= b.$$
Since $G[B_1 \cup B_2]$ is acyclic, $P_c$ has one end belonging to $B_3$. Hence, $G[X]$ and $G[Y]$ are acyclic, because $G[B_1 \cup B_2]$ is acyclic and $B_3$ is independent. Notice that $X$, $Y$ are disjoint, $abc \subset G[Y]$, ${(B_1  \backslash \{b\}) \cup B_2} \subset X$, $B_3 \cup \{b\} \subset Y$, and $X \cup Y = B_1 \cup B_2 \cup B_3 \cup V(P_a) \cup V(P_b)$. Therefore the conditions (a)--(c) yield.

If  $b \notin V(P_a) \cup V(P_c)$ and $P_a \neq P_c$, then we set
$$X =  B_1 \backslash \{b\} \cup B_2 \cup (V(P_a) \cap S_2 ) \cup (V(P_c) \cap S_2) ,$$
$$ Y = B_3 \cup \{b\} \cup  ( V(P_a) \cap S_3) \cup (V(P_c) \cap S_3)\quad \hbox{and} \quad v:= b.$$
If $b \in V(P_a) \cap V(P_c)$, then we set
$$X = B_1 \backslash \{b\} \cup B_2, \quad Y = B_3 \cup V(P_a) \cup V(P_c) \quad \hbox{and} \quad v:= b.$$
If $b \notin V(P_a)$ and $P_a = P_c$, then we set
$$X = B_1 \backslash \{b\} \cup B_2 \cup (V(P_a) \cap S_2), \quad Y = B_3 \cup \{b\} \cup (V(P_a) \cap S_3) \quad \hbox{and} \quad v:= b.$$

Case ($\beta, 2$). If $a,c \notin V(P_b)$, then we set
$$X = B_1 \cup B_2 \cup V(P_b - b) \quad \hbox{and} \quad Y = B_3 \cup  \{b\}.$$
If $a,c \in V(P_b)$, then  we set
$$X = B_1 \cup B_2  \quad \hbox{and} \quad  Y = B_3 \cup V(P_b)\backslash B_1.$$

Case ($\beta, 3$).
Let  $b \in V(P_a)$. Since  $G[B_1 \cup B_2]$ is acyclic, there is a vertex
$d \in V(P_a) \cap S_1$. Then we set
$$X = B_1 \backslash \{b\} \cup  B_2 \cup  \{d\}, \quad Y = B_3 \cup V(P_a)\backslash \{d\} \quad \hbox{and} \quad v:= b.$$
There are vertices $e,f \in B_2$ which are adjacent to every vertex of $P_a$. Since $G[B_1 \cup B_2]$ is acyclic, $ebf$ is the only $ef$-path consisting of vertices in $B_1 \cup B_2$. Hence, $G[X]$ and $G[Y]$ are acyclic, because $G[B_1 \cup B_2]$ is acyclic and $B_3$ is independent. Notice that $X$, $Y$ are disjoint, $abc \subset G[X]$, $({B_1 \backslash \{b\}) \cup B_2} \subset X$, $B_3 \cup \{b\} \subset Y$, and $X \cup Y = B_1 \cup B_2 \cup B_3 \cup V(P_b)$. Therefore the conditions (a)--(c) yield.

If $b \notin V(P_a)$, then we set
$$X = B_1 \backslash \{b\} \cup B_2 \cup (V(P_a) \cap S_2), \quad Y = B_3 \cup (V(P_a) \cap S_3) \cup \{b\} \quad \hbox{and} \quad v:= b.$$
Case ($\gamma$). For the cases ($\beta,1$)--($\beta,3$), Lemma~\ref{lemma2.3} follows from the assumption that the induced graph $G[B_1 \cup B_2]$ is acyclic. Hence, for the case ($\gamma$), Lemma~\ref{lemma2.3} yields.
\end{PrfFact}

\begin{theorem}\label{theorem2.1}
If induced graphs $G[B_1 \cup B_2]$ and $G[B_1 \cup B_3]$ are acyclic, then the following properties are satisfied:

(a) For every path $abc$ there is possible to partition the vertex set of $G$ into two subsets so that each induces a tree, and one of them contains the edge $ab$ and avoids the vertex $c$

(b)  For every path $abc$  with vertices $a$, $c$ of the same color there is possible to partition the vertex set of $G$ into two subsets so that each induces a tree, and one of them contains the path $abc$.
\end{theorem}
\begin{PrfFact}. Fix a path $abc$. If the vertices $a$, $c$ have different colors, the condition (a) follows from (b).
Therefore we assume that $a$, $c$ have the same color. If $G$ has at most two big vertices, then conditions (a) and (b) hold. Thus, assume that $G$ has at least three big vertices.

Suppose that $G$ is not $4$-connected, then there is a separating induced $3$-cycle, say $C$. Thus, $G-C$ consists of two components of $G$, say $C^+$ and $C^-$. Hence, by induction, for graphs $G_1 = G - C^{+}$ and $G_2 = G - C^{-}$ the conditions (a) and (b) follow. Let $d$ be a vertex of $C$ which has the same color as $a$ and $c$. If $ab$ is an edge in $G_1$, and $bc$ is an edge in $G_2$, then there exits a partitioning of the vertex set of $G_1$ ($G_2$) into two subsets $X_1$ and $Y_1$  ($X_2$ and $Y_2$) each of which induces a tree in $G_1$ ($G_2$)  such that the tree $G_{1}[X_1]$ ($G_{2}[X_2]$) contains the path $abd$ ($dbc$, respectively). If $abc$ is contained in $G_1$, then there exists a partitioning of $G_1$ ($G_2$) into two subsets $X_1$ and $Y_1$  ($X_2$ and $Y_2$) each of which induces a tree in $G_1$ ($G_2$) such that the tree $G_{1}[X_1]$ contains the path $abc$ and $X_1 \cap V(C) =  X_2 \cap V(C)$. In the both cases a partitioning $\{X_1 \cup X_2, Y_1 \cup Y_2\}$  of the vertex set of $G$ satisfies the condition~(b). By analogy, there exists a partitioning of the vertex set of $G$ satisfying the condition (a).

If $G$ is $4$-connected with at least three big vertices, the conditions (a) and (b) follow from Lemmas~\ref{lemma2.1}--\ref{lemma2.2}--~\ref{lemma2.3}.
\end{PrfFact}

\end{document}